\documentclass[reqno,centertags, 12pt]{amsart}
\usepackage{amsmath,amsthm,amscd,amssymb}
\usepackage{latexsym}
\newcommand{\bbR}{{\mathbb{R}}}

\newcommand{\bbZ}{{\mathbb{Z}}}

\newcommand{\calH}{{\mathcal H}}

\newcommand{\lb}{\label}
\newcommand{\f}{\frac}

\newcommand{\ti}{\tilde  }
\newcommand{\wti}{\widetilde  }

\newcommand{\dist}{\text{\rm{dist}}}
\newcommand{\loc}{\text{\rm{loc}}}

\newcommand{\spec}{\text{\rm{spec}}}

\newcommand{\bi}{\bibitem}

\newcommand{\beq}{\begin{equation}}
\newcommand{\eeq}{\end{equation}}
\newcommand{\ba}{\begin{align}}
\newcommand{\ea}{\end{align}}
\newcommand{\veps}{\varepsilon}

\newcounter{smalllist}

\allowdisplaybreaks
\numberwithin{equation}{section}

\newtheorem{theorem}{Theorem}[section]

\newtheorem*{p2.1}{Proposition 2.1}

\newtheorem{lemma}[theorem]{Lemma}
\newtheorem{corollary}[theorem]{Corollary}
\theoremstyle{definition}
\newtheorem{example}[theorem]{Example}

\theoremstyle{remark}
\newtheorem*{remark}{Remark}
\newtheorem*{remarks}{Remarks}

\newcommand{\abs}[1]{\lvert#1\rvert}
\newcommand{\jap}[1]{\langle #1 \rangle}

\begin{document}

\title[Eigenvalue Bounds]{Eigenvalue Bounds for
Perturbations of Schr\"odinger Operators and Jacobi Matrices With Regular Ground States}
\author[R.~Frank, B.~Simon, and T.~Weidl]{Rupert L.~Frank$^1$, Barry Simon$^2$, and
Timo Weidl$^3$}

\thanks{$^1$ Department of Mathematics, Royal Institute of Technology, 100~44 Stockholm,
Sweden. E-mail: rupert@math.kth.se}

\thanks{$^2$ Mathematics 253-37, California Institute of Technology, Pasadena, CA 91125, USA.
E-mail: bsimon@caltech.edu. Supported in part by NSF Grant DMS-0140592 and
U.S.--Israel Binational Science Foundation (BSF) Grant No.\ 2002068}

\thanks{$^3$ Stuttgart University, Department of Mathematics and Physics, Pfaffenwaldring~57,
70569~Stuttgart, Germany. Email: Timo.Weidl@mathematik.uni-stuttgart.de. Supported in part by
DFG grant WE-1964 2/1}

\date{June 27, 2007}
\keywords{Schr\"odinger operators, Jacobi matrices, Lieb--Thirring bounds}
\subjclass[2000]{81Q05, 42C05, 34L15}

\begin{abstract} We prove general comparison theorems for eigenvalues of perturbed
Schr\"odinger operators that allow proof of Lieb--Thirring bounds for suitable non-free
Schr\"odinger operators and Jacobi matrices.
\end{abstract}

\maketitle

\section{Introduction} \lb{s1}

Consider a Schr\"odinger operator,
\begin{equation} \lb{1.1}
H_0 = -\Delta + V_0
\end{equation}
on $L^2(\bbR^\nu)$ which we suppose obeys
\begin{equation} \lb{1.2}
\inf\spec(-\Delta + V_0)=0
\end{equation}
(By subtracting a constant, we can always arrange this, and by assuming this, the notation
simplifies.) We are interested in controlling the negative eigenvalues of
\begin{equation} \lb{1.3}
H=H_0+V
\end{equation}
We let
\begin{equation} \lb{1.4}
E_1 (V_0;V)\leq E_2(V_0;V) \leq \cdots\leq E_n(V_0;V)\leq \cdots
\end{equation}
be either the negative eigenvalues or $0$, that is,
\begin{equation} \lb{1.5}
E_j(V_0;V)=\min(0;\inf \{\lambda,\dim P_{(-\infty,\lambda]}(H)\geq j\})
\end{equation}

We say that $H_0$ has a {\it regular ground state\/} if and only if there exists a function, $u_0$,
on $\bbR^\nu$ obeying
\begin{gather}
(-\Delta +V_0) u_0=0 \lb{1.6} \\
0 < c_1 \leq u(x) \leq c_2 < \infty \lb{1.7}
\end{gather}
for some $c_1,c_2$. We take $c_1=\inf u$, $c_2 =\sup u$, and let
\begin{equation} \lb{1.8}
\beta (V_0)=\biggl( \f{c_2}{c_1}\biggr)^2
\end{equation}

Our main result is that any bound on the number or sums of eigenvalues for the operator $-\Delta+V$ can be carried
(with a change in the coupling constant) to the operator $H_0+V$\!. This is based on the following observation:

\begin{theorem}\lb{T1.1} For any $V\leq 0$ and $V_0$ obeying \eqref{1.2} with regular ground state, we have for
all $j$,
\begin{equation} \lb{1.9}
\abs{E_j (0;\beta^{-1}V)} \leq \abs{E_j (V_0;V)} \leq \abs{E_j (0;\beta V)}
\end{equation}
\end{theorem}

This result is remarkable for its generality and also for the simplicity of its proof. We will see
in Section~\ref{s2} that it can be used to compare not only $V_0$ and $0$ but two arbitrary $V_0$'s
with relatively bounded ground states.

Of course, \eqref{1.9} immediately implies bounds on moments of bound states:
\begin{equation} \lb{1.10}
S_\gamma (V_0;V) =\sum_{j=1}^\infty \, \abs{E_j (V_0;V)}^\gamma
\end{equation}
where we look at this only for $\gamma\geq 0$ and interpret $0^0 =0$ so $S_0$ is the number of strictly
negative eigenvalues. Clearly, Theorem~\ref{T1.1} implies:

\begin{corollary}\lb{C1.2} For any $\gamma >0$ and $\beta=\beta(V_0)$ given by \eqref{1.8},
\begin{equation} \lb{1.11}
S_\gamma (0;\beta^{-1} V) \leq S_\gamma (V_0;V) \leq S_\gamma (0;\beta V)
\end{equation}
\end{corollary}

The standard Lieb--Thirring inequalities (reviewed in \cite{Hun2007,LW}) assert
\begin{equation} \lb{1.12}
S_\gamma (0;\beta V)\leq L_{\gamma,\nu}\int \abs{V(x)}^{\gamma + \f{\nu}{2}}\, d^\nu x
\end{equation}
for $\gamma\geq\f12$ in $\nu=1$, $\gamma >0$ in $\nu=2$, and $\gamma \geq 0$ in $\nu\geq 3$. In some
cases, the optimal constants are known, for example, $L_{\f{1}{2},1}=\f12$. (These yield good
constants in our perturbed estimates but we do not claim optimal constants for our situation!)
Clearly, \eqref{1.11} implies:

\begin{corollary}\lb{C1.3}
\begin{equation} \lb{1.13}
S_\gamma (V_0;V)\leq L_{\gamma,\nu} \beta^{\gamma +\f{\nu}{2}}
\int \abs{V(x)}^{\gamma+ \f{\nu}{2}}\, d^\nu x
\end{equation}
and, in particular,
\begin{equation} \lb{1.14}
\sum_{j=1}^\infty\, \abs{E_j (V_0;V)}^{\f{1}{2}} \leq \tfrac12\, \beta \int \abs{V(x)}\, dx
\end{equation}
in $\nu=1$ dimension.
\end{corollary}

One can also obtain logarithmic estimates as in \cite{KVW} and Hardy--Lieb--Thirring bounds
as in \cite{EF}. Since it is known \cite{ZF,Sch,LW} that for $\nu=1$ and $V\leq 0$,
\begin{equation} \lb{1.15}
S_{\f{1}{2}} (0;V)\geq \tfrac14\, \int \abs{V(x)}\, dx
\end{equation}
we conclude that
\begin{equation} \lb{1.16}
S_{\f{1}{2}} (V_0;V) \geq \f{1}{4\beta}\, \int \abs{V(x)}\, dx
\end{equation}

These results are of interest because there are many cases which are known to have regular ground
states.

\begin{example}\lb{E1.4} If $V_0$ is periodic, then there is a positive periodic ground state.
If $V_0$ is locally $L^{\f{\nu}{2}}$ (if $\nu\geq 3$, locally $L^1$ if $\nu=1$, and locally $L^p$
with $p>1$ if $\nu=2$), then it is known that eigenfunctions are continuous (see \cite{Sxxi})
and thus, $H_0$ has a regular ground state.
\qed
\end{example}

\begin{example} \lb{E1.5} We will discuss Jacobi matrices in Sections~\ref{s3} and \ref{s4}. It is
known (see \cite{SY,PY,CSZ}) that elements in the isospectral torus of finite gap Jacobi matrices
have regular ground states.
\qed
\end{example}

\begin{example}\lb{E1.6} If $u$ is any function obeying \eqref{1.7}, then $V_0= \f{\Delta u}{u}$ has
a regular ground state.
\qed
\end{example}

In Section~\ref{s2}, we will review the ground state representation and prove a stronger theorem
than Theorem~\ref{T1.1}. As hinted, it is the ground state representation that is critical. In
this regard, we should emphasize that the variational argument we use in Section~\ref{s2} has appeared
earlier in work of the Birman school---we would mention, in particular, Lemma~6.1 of Birman, Laptev,
and Suslina \cite{BLS}, although it may have appeared earlier in their work. Our novelty here is
the wide applicability, the use in CLR and Lieb--Thirring bounds, and the applicability to the
discrete case and Szeg\H{o} estimates.

As we will explain in Section~\ref{s4}, an initial motivation for this work was critical Lieb--Thirring
bounds for finite gap almost periodic Jacobi matrices in connection with Szeg\H{o}'s theorem for
such situations. Ground state representations do not seem to be in the literature for Jacobi
matrices, so we do this first in Section~\ref{s3}, and then prove an analog of Theorem~\ref{T1.1}
for Jacobi matrices in Section~\ref{s4}. Section~\ref{s5} discusses some other cases.

\smallskip
It is a pleasure to thank Fritz Gesztesy, Yehuda Pinchover, Robert Seiringer, and
Simone Warzel for useful comments, and Michael Aizenman for being a sensitive editor.

\section{Comparison for Schr\"odinger Operators} \lb{s2}

Fundamental to our results is the ground state representation that if \eqref{1.6} holds for $u_0$,
continuous and strictly positive on $\bbR^\nu$, then
\begin{equation} \lb{2.1}
\jap{gu_0, H_0 gu_0} = \int \abs{\nabla g}^2 u_0^2\, d^\nu x
\end{equation}

Ground state representations go back to Jacobi \cite{Jac}. For Schr\"odinger operators, it appears at least
as far back as Birman \cite{Bir} and it was used extensively in constructive quantum field theory
(especially by Segal, Nelson, Gross, and Glimm--Jaffe; see Glimm--Jaffe \cite{GJ}). As a basis for
comparison theorems, it was used by Kirsch--Simon \cite{KS} and, as noted above, in a similar context
by Birman--Laptev--Suslina \cite{BLS}.

We will be cavalier about technical assumptions needed for \eqref{2.1}. From one point of view, we can
use \eqref{2.1} as a definition of $H_0$! Namely, the right side of \eqref{2.1} defined for $g$'s
with distributional derivative making the right side finite is easily seen to be a closed quadratic
form on $\calH_{u_0}\equiv L^2 (\bbR^\nu, u_0^2 d^\nu x)$ defining a positive selfadjoint operator
$\wti H_0$ on $\calH_{u_0}$. The unitary operator $W\colon L^2 (\bbR^\nu, d^\nu x)\to\calH_{u_0}$
by $Wg=u_0^{-1} g$ lets us define $H_0=W^{-1} \wti H_0 W$ and our results hold for perturbations
of that.

It is not hard to prove that if $V_0=V_{0,+} + V_{0,-}$ with $V_{0,+}\in L_\loc^1 (\bbR^\nu, d^\nu x)$
and $V_{0,-}\in K_\nu$, the Kato class, then the selfadjoint operator $H_0$ defined as the form
closure of $-\Delta +V_0$ on $C_0^\infty$ obeys \eqref{2.1} if $u_0$ is a positive distributional
solution of \eqref{1.6}.

Notice that we do not need \eqref{1.2}, but only $\inf\spec (H_0)\geq 0$ for this to work, and
Theorem~\ref{T2.1} below holds in that case (although $\inf\spec(H_0)=\inf\spec (H_1)$) under the
hypothesis of the theorem. Here is our main result:

\begin{theorem}\lb{T2.1} Let $H_0,H_1$ have the form \eqref{2.1} for positive continuous functions
$u_0,u_1$. Suppose
\begin{equation} \lb{2.2}
0 < \inf\biggl( \f{u_0}{u_1}\biggr) \leq  \sup \biggl( \f{u_0}{u_1}\biggr) <\infty
\end{equation}
and let
\begin{equation} \lb{2.3}
\beta \equiv \biggl[ \f{\sup (\f{u_0}{u_1})}{\inf (\f{u_0}{u_1})}\biggr]^2
\end{equation}
For any $V\leq 0$, let $E_j$ be given by \eqref{1.5}. Then
\begin{equation} \lb{2.4}
\abs{E_j (V_0;V)}\leq \abs{E_j (V_1;\beta V)}
\end{equation}
\end{theorem}

\begin{remark} By interchanging $V_0$ and $V_1$ and replacing $V$ by $\beta^{-1} V$\!, we get the
complementary inequality
\begin{equation} \lb{2.5}
\abs{E_j(V_1;\beta^{-1} V)} \leq \abs{E_j (V_0;V)}
\end{equation}
\end{remark}

\begin{lemma}\lb{L2.1} Let $V\leq 0$. Let $\tau >0$. If for some $g$,
\begin{equation} \lb{2.6}
\jap{gu_0, (H_0+V)gu_0}\leq -\tau \jap{gu_0,gu_0}
\end{equation}
then
\begin{equation} \lb{2.7}
\jap{gu_1, (H_1+\beta V) gu_1} \leq -\tau \jap{gu_1, gu_1}
\end{equation}
\end{lemma}

\begin{proof} Let
\begin{align}
\beta_+ &= \sup\biggl(\f{u_0}{u_1}\biggr)^2 = \biggl[\inf \biggl(\f{u_1}{u_0}\biggr)\biggr]^{-2} \lb{2.8} \\
\beta_- &= \inf\biggl(\f{u_0}{u_1}\biggr)^2 = \biggl[\sup\biggl(\f{u_1}{u_0}\biggr)\biggr]^{-2} \lb{2.9}
\end{align}
so
\begin{equation} \lb{2.10}
\beta_+ = \beta\beta_- \Rightarrow \beta\beta_+^{-1} =\beta_-^{-1}
\end{equation}

Since $-Vg^2 \geq 0$ and $\abs{\nabla g}^2\geq 0$, we have
\begin{align}
\int(\nabla g)^2 u_1^2\, dx &\leq \beta_-^{-1} \int (\nabla g)^2 u_0^2\, dx \lb{2.11} \\
-\int Vg^2 u_1^2\, dx &\geq -\beta_+^{-1} \int Vg^2 u_0^2\, dx \lb{2.12}
\end{align}
so, by \eqref{2.10} and \eqref{2.1},
\begin{align}
\jap{gu_1, (H_1+\beta V)gu_1} &\leq \beta_-^{-1} \jap{gu_0, (H_0 +V)gu_1} \notag \\
&\leq -\tau\beta_-^{-1} \jap{gu_0, gu_0} \lb{2.13}
\end{align}

But
\begin{equation} \lb{2.14}
\jap{gu_0, gu_0}\geq \beta_- \jap{gu_1, gu_1}
\end{equation}
and $\tau >0$, so
\[
\text{RHS of \eqref{2.13}} \leq -\tau \jap{gu_1,gu_1}
\]
proving \eqref{2.7}.
\end{proof}

\begin{proof}[Proof of Theorem~\ref{T2.1}] If $E_j(V_0;V)=0$, there is nothing to prove. If $\tau\equiv
\abs{E_j (V_0;V)} >0$, there is a space, $\calH_j$, of dimension at least $j$ so
\begin{equation} \lb{2.15}
\jap{\psi, (H_0+V)\psi}\leq -\tau\jap{\psi,\psi}
\end{equation}
for $\psi\in\calH_j$. By the lemma, if $\varphi =\f{u_1}{u_0}\psi$, we have
\begin{equation} \lb{2.16}
\jap{\varphi, (H_1 +\beta V)\varphi} \leq -\tau\jap{\varphi,\varphi}
\end{equation}

Thus, there is a space of dimension at least $j$ where \eqref{2.16} holds. By the min-max
principle (see \cite{RS4}),
\begin{equation} \lb{2.17}
E_j (V_1;\beta V)\leq -\tau
\end{equation}
which is \eqref{2.4}.
\end{proof}

\section{Ground State Representation for Jacobi Matrices} \lb{s3}

While we are interested mainly in semi-infinite one dimension Jacobi matrices,
that is, tridiagonal semi-infinite matrices, we can consider the higher-dimensional
case as well, so we will. So far as we know, there is no prior literature on the
ground state representation for discrete operators, so we start with that in
this section.

In $\bbZ^\nu$, we let $\delta_j$, $j=1,\dots,\nu$, be the $\nu$ component vectors
with $1$ in the $j$th place and $0$ elsewhere. So $k\pm\delta_j$ are the $2\nu$ neighbors
of $k\in\bbZ^\nu$. A Jacobi operator is parametrized by a symmetric $a_{j\ell} >0$ for
each $j$, $\ell\in\bbZ_\nu$ with $\abs{j-\ell}=1$ and a real number $b_k$ for each
$k\in\bbZ^\nu$. We will suppose
\[
\sup_k\, \abs{b_k} + \sup_\ell\, \abs{a_{k\ell}} <\infty
\]

The Jacobi operator associated with these parameters is the operator $J$ on $\ell^2
(\bbZ^\nu)$ with
\begin{equation} \lb{3.1}
(J\varphi)_\ell = \sum_{\pm, j=1,\dots,\nu} a_{\ell\, \ell\pm\delta_j} \varphi_{\ell\pm\delta_j}
+b_\ell \varphi_\ell
\end{equation}
We will use $J(\{a_{\ell m}\}, \{b_\ell\})$ if we want to make the dependence on $a$ and $b$ explicit.

\begin{lemma}\lb{L3.1} Let $f$ be a bounded real-valued function on $\bbZ^\nu$ and $M_f$ the diagonal
matrix on $\ell^2 (\bbZ^\nu)$ which is multiplication by $f$. Then
\begin{align}
[M_f, J(\{a_{\ell m}\}, \{b_\ell\})] &= J(\{a_{\ell m}(f_\ell-f_m), b_\ell\equiv 0\}) \lb{3.2} \\
[M_f, [M_f, J(\{a_{\ell m}\}, \{b_\ell\})]] &= J(\{a_{\ell m} (f_\ell-f_m)^2, b_\ell\equiv 0\}) \lb{3.3}
\end{align}
\end{lemma}

\begin{proof} \eqref{3.2} is an elementary calculation and it implies \eqref{3.3}.
\end{proof}

\begin{theorem}\lb{T3.2} Let $J$ be a Jacobi operator on $\ell^2 (\bbZ^\nu)$ with parameters
$\{a_{\ell m}\}, \{b_\ell\}$. Suppose $u$ is a positive ``sequence" parametrized by $\bbZ^\nu$
so that
\begin{equation} \lb{3.4}
\sum_{\abs{m-\ell}=1} a_{\ell m} u_m + b_\ell u_\ell =0
\end{equation}
for all $\ell\in\bbZ^\nu$. Then for any $f$ with $fu\in \ell^2 (\bbZ^\nu)$, we have
\begin{equation} \lb{3.5}
\jap{fu, (-J)fu} =\sum_{\substack{m,\ell \\ \abs{m-\ell}=1}} a_{\ell m} u_\ell u_m (f_\ell -f_m)^2
\end{equation}
\end{theorem}

\begin{remark} In particular, this shows $-J\geq 0$.
\end{remark}

\begin{proof} It suffices to prove \eqref{3.5} for $f$ of finite support and take limits. For the
left side converges since $-J$ is a bounded operator, and by positivity, the right side converges
(a priori perhaps to $\infty$, but by the equality to a finite limit; it is only here that positivity
of $u$ is used).

For finite sequences, $f$, use \eqref{3.3}, taking expectations in a vector $\ti u$ which equals $u$
on $\{\ell\mid f(m)$ is non-zero for some $m$ with $\abs{m-\ell}\leq 1\}$. Then $M_f J\ti u=0$, so
$\jap{\ti u, [M_f, [M_f J]]\ti u} =-2 \jap{fu,Jfu}$.
\end{proof}

\section{Comparison for Jacobi Operators} \lb{s4}

In this section, we will prove an analog of Theorem~\ref{T1.1} for Jacobi operators. One difference
that we have to expect is that of $\sup\spec(J)=0$; the same is true of $\lambda J$ for any $\lambda
>0$, but the ground states are the same. Thus, comparison of the two $J$'s cannot involve only
the ground state ratio but also a setting of scales which will enter as a ratio of $a$'s. In the
Schr\"odinger case, the scale is set by the $-\Delta$ rather than $-\lambda\Delta$.

For notation, we let $E_j (\{a_{\ell m},b_\ell\})$ be the max of zero and the $j$th eigenvalue of
$J(\{a_{\ell m},b_\ell\})$ counting from the top. Here is our main result:

\begin{theorem} \lb{T4.1} Let $\{a_{\ell m}^{(0)}, b_\ell^{(0)}\}$ and $\{a_{\ell m}^{(1)}, b_\ell^{(1)}\}$
be two sets of bounded Jacobi parameters with positive sequences $u^{(0)}, u^{(1)}$ obeying \eqref{3.4}
for $\{a^{(0)}, b^{(0)}\}, \{a^{(1)}, b^{(1)}\}$, respectively. Let
\begin{align}
\beta_+ &= \sup_\ell \biggl( \f{u_\ell^{(0)}}{u_\ell^{(1)}}\biggr)^2 \lb{4.1} \\
\beta_- &= \inf_\ell \biggl( \f{u_\ell^{(0)}}{u_\ell^{(1)}}\biggr)^2 \lb{4.2} \\
\gamma_- &= \inf \biggl( \f{a_{j\ell}^{(0)} u_j^{(0)} u_\ell^{(0)}}
{a_{j\ell}^{(1)} u_j^{(1)} u_\ell^{(1)}}\biggr) \lb{4.3}
\end{align}
Suppose $\beta_+ <\infty$ and $\beta_- >0$. Then for perturbations $\{\delta a_{\ell m}, \delta b_\ell\}$
with $a_{\ell m}^{(0)} + \delta a_{\ell m}>0$,
\begin{equation} \lb{4.4}
E_j (\{a_{\ell m}^{(0)} + \delta a_{\ell m}, b_\ell^{(0)} + \delta b_\ell\}) \leq
E_j \biggl( \biggr\{\eta a_{\ell m}^{(1)}, \eta b_\ell^{(1)} + \beta \biggl[
\abs{\delta b_\ell} + \sum_{\abs{m-\ell} =1} \, \abs{\delta a_{\ell m}}\biggr] \biggr\}\biggr)
\end{equation}
where
\begin{equation} \lb{4.5}
\eta = \f{\gamma_-}{\beta_-} \qquad
\beta = \f{\beta_+}{\beta_-}
\end{equation}
\end{theorem}

\begin{remarks} 1. We only have a one-sided inequality as we would in the Schr\"odinger case
if we did not demand $V\leq 0$. Since $\delta a$ terms are never of a definite sign, we cannot
have them in a two-sided comparison. But there is clearly a two-sided comparison if $\delta a=0$
and $\delta b>0$.

\smallskip
2. Note that rescaling $u^{(0)}$ or $u^{(1)}$ which changes $\beta_+,\beta_-,\gamma_-$ does not
change $\eta$ or $\beta$. Similarly, $\eta$ scales properly under changes of the scale of $a$.
\end{remarks}

Following Hundertmark--Simon \cite{HunS}, we begin with a reduction to the case $\delta a=0$,
$\delta b\geq 0$:

\begin{lemma}\lb{L4.2} One has
\begin{equation} \lb{4.6}
E_j (\{a_{\ell m}^{(0)} + \delta a_{\ell m}, b_\ell^{(0)} +\delta b_\ell\}) \leq
E_j \biggl(\biggl\{a_{\ell m}^{(0)}, b_\ell^{(0)} + \abs{\delta b_\ell} +
\sum_{\abs{m-\ell}=1} \abs{\delta a_{\ell m}}\biggr\}\biggr)
\end{equation}
Thus, it suffices to prove \eqref{4.4} when $\delta a_{\ell m}=0$ and $\delta b_\ell \geq 0$.
\end{lemma}

As in \cite{HunS}, this follows from
\[
\begin{pmatrix} 0 & a \\ a & 0 \end{pmatrix} \leq
\begin{pmatrix} a & 0 \\ 0 & a \end{pmatrix}
\]
and $b\leq \abs{b}$.

\begin{proof}[Proof of Theorem~\ref{T4.1}] The proof is identical to the proof in Section~\ref{s2},
the sole change being that \eqref{2.11} needs to be replaced by
\begin{equation} \lb{4.7}
\sum_{\substack{\ell, m \\ \abs{\ell-m}=1}} a_{\ell m}^{(1)} u_\ell^{(1)} u_m^{(1)} \abs{f_\ell -f_m}^2
\leq \gamma_-^{-1} \sum_{\substack{\ell, m \\ \abs{\ell-m} =1}} a_{\ell m}^{(0)} u_\ell ^{(0)}
u_m^{(0)} \abs{f_\ell -f_m}^2
\end{equation}
So instead of a $\f{\beta_-}{\beta_-}=1$, we get $\eta = \f{\gamma_-}{\beta_-}$.
\end{proof}

Let $W\colon \ell^2 (\bbZ^\nu)\to \ell^2 (\bbZ^\nu)$ by
\begin{equation} \lb{4.8}
(Wf)_n =(-1)^{\abs{n}} f_n
\end{equation}
with $\abs{n}=\abs{n_1} + \cdots + \abs{n_\nu}$. Then
\begin{equation} \lb{4.9}
W\!J(\{a_{\ell m}, b_\ell\}) W^{-1} = -J(\{a_{\ell m}, -b_\ell\})
\end{equation}
This allows one to control eigenvalues below $\inf\spec(J)$ in the same way. We define $E_j^-$ to be
the min of $0$ and the $j$th eigenvalue of $J$ counting from the bottom. We call a sequence $u$
on $\bbZ^\nu$ alternating positive if $(-1)^{\abs{m}} u_m>0$ for all $m$. The result is:

\begin{theorem}\lb{T4.3} Let $\{a_{\ell m}^{(0)}, b_\ell^{(0)}\}$ and $\{a_{\ell m}^{(1)},
b_\ell^{(1)}\}$ be two sets of bounded Jacobi parameters with alternating positive sequences $u_m^{(0)},
u_m^{(1)}$ obeying \eqref{3.4} for $\{a^{(0)},b^{(0)}\}$, $\{a^{(1)}, b^{(1)}\}$, respectively.
Define $\beta_+, \beta_-, \gamma_-, \eta, \beta$ as in \eqref{4.1}--\eqref{4.3} and \eqref{4.5}.
Then \eqref{4.4} holds with $E_j$ replaced by $\abs{E_j^-}$.
\end{theorem}

\begin{corollary} Let $\{a_{\ell m}^{(0)}, b_\ell^{(0)}\}$ be a one dimension periodic set of
Jacobi parameters or the almost periodic parameters associated with a finite gap spectrum
{\rm{(}}see \cite{SY}{\rm{)}}. Let $J_0$ be the associated half-line Jacobi matrix. Let $J$ be
the Jacobi matrix associated with $\{a_{\ell m}^{(0)} + \delta a_{\ell m}, b_\ell^{(0)} +
\delta b_\ell\}$ where
\begin{equation} \lb{4.10}
\sum_{\substack{ \abs{\ell-m}=1 \\ \ell=1,2,3, \dots }} \, \abs{\delta a_{\ell m}} +
\sum_{\ell=1}^\infty \, \abs{\delta b_\ell}<\infty
\end{equation}
Let $E_1^- < E_2^- < \cdots < \inf\spec (J_0) < \sup\spec(J_0) < \cdots < E_2^+ < E_1^+$
be the eigenvalues of $J$ outside the convex hull of $\spec(J_0)$. Then
\begin{equation} \lb{4.11}
\sum_{k,\pm} \dist (E_k^\pm, \sigma(J_0))^{\f{1}{2}} \leq C \sum (\abs{\delta a_{\ell m}} +
\abs{\delta b_\ell})
\end{equation}
for a constant $C$ depending only on $J_0$.
\end{corollary}

\begin{proof} We only need that $J_0$ and $W\!J_0W^{-1}$ have regular ground states. This follows from Floquet
theory for the periodic case and from the detailed analysis of Jost solutions for the almost periodic case;
see \cite{SY,PY,CSZ}. Then compare to the free $J_0$ ($a_{\ell m}=1$ if $\abs{\ell-m}=1$, $b_\ell=0$)
and use the bound of \cite{HunS}.
\end{proof}

For the periodic case, this is proven by Damanik, Killip, and Simon \cite{DKS06}, who also prove
this where the sum in \eqref{4.11} is over all eigenvalues including the ones in gaps.

It remains an interesting question relevant to the study of the Szeg\H{o} condition (see
\cite{HSppt,CSZ}) to get the bound in the almost periodic case. In \cite{HSppt},
Hundertmark and Simon prove the weaker bounds where the $\f12$ power is replaced by any
$p>\f12$ or where $p=\f12$ but there is an $\ell^\veps$ added to the sum.

\section{Some Final Remarks} \lb{s5}

All we needed for our arguments is some kind of ground state representation. That means we can
replace $-\Delta$ by
\begin{equation} \lb{5.1}
L_0 =-\sum_{j,k=1}^\nu \partial_j A_{jk}(x)\, \partial_k
\end{equation}
with $\{A_{jk}(x)\}_{1\leq j,k\leq \nu}$ a strictly positive matrix. Hence, if $(L_0+V)u_0=0$, then
\[
\jap{fu_0, (L_0+V) fu_0} =\int \jap{\nabla f, A\nabla f} \abs{u_0(x)}^2\, d^\nu x
\]
and we can still compare to $-\Delta$ although $\sup\|A\|$ and $\inf \|A\|$ will enter.

We can compare magnetic field operators where the magnetic field is fixed but $V_0,V$ vary. 
For if $H_0u\equiv (-(\nabla -ia)^2 +V_0)u_0=0$, then, since $[f^*, [f,H_0]] = 
\abs{(\nabla -ia)f}^2$, we have that
\[
\jap{fu_0, H_0fu_0} =\int \abs{(\nabla-ia)f}^2 \abs{u_0}^2 \, d^\nu x
\]

For the discrete case, the key was not tridiagonal matrices, but ones with a ground state representation.
For example, non-negative off-diagonal will do.

Using the representation of Frank, Lieb, and Seiringer \cite{FLS}, one can treat some perturbations of
$(-\Delta)^\alpha$, $0<\alpha <1$.


\end{document}